\documentclass[11pt]{article}
\usepackage{amsmath}
\usepackage{mathrsfs}
\usepackage{graphicx}
\def\bel{\begin{equation}\label}
\def\eeq{\end{equation}}
\newtheorem{Definition}{Definition}[section]
\newtheorem{Theorem}{Theorem}[section]
\newtheorem{Remark}{Remark}[section]
\newtheorem{Lemma}{Lemma}[section]
\newtheorem{Proposition}{Proposition}[section]

\newtheorem{Example}{Example}[section]

\title{Asymptotic controllability  and optimal control}
\author{M. Motta  \& F. Rampazzo\\
Dipartimento di Matematica\\
Via Trieste, 63 - 35121 Padova, Italy\\
Telefax (39)(049) 8271428\\
e-mail: motta@math.unipd.it \,\,\,\,rampazzo@math.unipd.it}

\def\fudge{\mathchoice{}{}{\mkern.5mu}{\mkern.8mu}}
\def\bbc#1#2{{\rm \mkern#2mu\vbar\mkern-#2mu#1}}
\def\bbb#1{{\rm I\mkern-3.5mu #1}}
\def\bba#1#2{{\rm #1\mkern-#2mu\fudge #1}}
\def\bb#1{{\count4=`#1 \advance\count4by-64 \ifcase\count4\or\bba A{11.5}\or
\bbb B\or\bbc C{5}\or\bbb D\or\bbb E\or\bbb F \or\bbc G{5}\or\bbb H\or
\bbb I\or\bbc J{3}\or\bbb K\or\bbb L \or\bbb M\or\bbb N\or\bbc O{5} \or
\bbb P\or\bbc Q{5}\or\bbb R\or\bbc S{4.2}\or\bba T{10.5}\or\bbc U{5}\or%
\bbb P\or\bbc Q{5}\or\bbb R\or\bba S{8}\or\bba T{10.5}\or\bbc U{5}\or
\bba V{12}\or\bba W{16.5}\or\bba X{11}\or\bba Y{11.7}\or\bba Z{7.5}\fi}}

\def \R {{\bb R}}
\def \C {{\bf C}}
\def \tC {{\hat{\bf C}}}
\def \N {{\bb N}}

\def \V{{{\mathcal V}}}

\def\rr{I\!\!R}
\def\ds{\displaystyle}
\def \vv{\vskip 0.5 truecm}

\begin{document}
\maketitle
\begin{abstract}
{We consider a control problem where the state must  approach asymptotically a target $\C$ while paying an integral cost with a {\em non-negative} Lagrangian $l$. The dynamics $f$ is just continuous, and  no assumptions are made on the zero level set of the Lagrangian $l$.  Through an inequality  involving a positive number $\bar p_0$ and a {\it Minimum Restraint Function $U=U(x)$} --a special type of  Control Lyapunov Function--   we provide a   condition implying   that  {\bf (i)}  the system is    asymptotically controllable, and {\bf (ii)}  the value function is bounded by $U/\bar p_0$. The result has significant consequences for the uniqueness issue of the corresponding Hamilton-Jacobi equation. Furthermore it may be regarded as a first step in the direction of a feedback construction.}
\end{abstract}

\footnotetext  {$({\bf *})$ This research is partially supported by the Marie Curie ITN SADCO, FP7-PEOPLE-2010-ITN n. 264735-SADCO, by the MIUR grant PRIN 2009
"Metodi di viscosit\`a, geometrici e di controllo per modelli diffusivi non lineari" (2009KNZ5FK), and by the Fondazione CaRiPaRo Project "Nonlinear Partial Differential Equations: models, analysis, and control-theoretic problems".}

\section{Introduction}

 Let ${\bf C}\subset\R^n$ be a closed subset, which  will be called the {\it target}, and let $\C^c$ denote its complement. We consider the value function
\bel{P} \begin{array}{l}
{{\mathcal V}}(x)\doteq{\inf} \, {\mathcal J}_{z,\alpha},\,\,
\\\,\\
 \ds{\mathcal J}_{z,\alpha} \doteq \int_0^{{T}_{z}}  l(z(t),\alpha(t))\,dt
 \end{array}
 \eeq
for trajectory-control pairs $(z,\alpha):[0,{T}_{z}[:\to \C^c\times A$ subject to
\bel{ode}\begin{array}{l}\dot z(t)=f(z(t),\alpha(t)) \qquad z(0)=x\\\,\\ \lim_{t\to {T}_{z}^-}{\bf d}(z(t),{\bf C})=0,\end{array}\eeq
where ${\bf d}$ denotes the Euclidean distance.

The crucial assumption of the paper will be the sign condition
   \bel{sign} l(x,a)\geq 0 \qquad \forall (x,a)\in \C^c\times A.\eeq

   In its stronger form, the main result of the paper (Theorem \ref{Th3.1gen}) reads as follows:
\vskip0.3truecm
    {\it Let $\bar p_0$ be a positive real number and let $U=U(x)$ be a proper, positive definite, semiconcave function such that, for every $x\in\C^c$,  one has
   \begin{equation}\label{C11}
H(x,\bar p_0,D^* U(x))<0.  \,\,\,\, \footnotemark
\end{equation}
\footnotetext{The notation  $H(x,\bar p_0,D^* U(x))<0$  means that $H(x,\bar p_0,p)<0$  $\forall p\in D^* U(x)$.}
 Then the system is globally asymptotically controllable and  the value function $\V$ verifies the inequality $$\V(x)\leq U(x)/\bar p_0,$$ for all $x\in\C^c$.}
 \vskip0.4truecm
 In (\ref{C11}), $H$ denotes the natural (minimized)  Hamiltonian of the system, namely\footnotemark
   \bel{Hdef}H(x,p_0,p)\doteq \inf_{a\in A}\left\langle (p_0, p) \,,\,\Big(l(x,a),f(x,a)\Big)\right\rangle,\,\,\,\,\,
\eeq

\footnotetext{The fact that the state variable is $n$-dimensional while the adjoint variable $(p_0,p)$ is $(n+1)$-dimensional is due to the presence of a  {\it hidden} state variable $z_0$, namely the one verifying the differential equation $\dot z_0(t) = l(z(t),\alpha(t))$.}
 while  $D^*$ is  the generalized differential operator called  {\it limiting gradient} (see Def. \ref{D*}).

   A function $U$   as above  is here  called a  {\it Minimum Restraint Function} (Def. \ref{CMRL}).

Let us make clear that condition (\ref{C11}) {\it is not} a mere application of the usual (first order) asymptotic global controllability condition  to the enlarged dynamics obtained by adding the equation $\dot z_0=l(z,\alpha)$, with  the enlarged target $[0,+\infty[\times\C$. Actually, the  known conditions on the existence of a Control Lyapunov Function to characterize global asymptotic controllability would  provide no information on the value of the minimum \footnote{ Notice also that $W(x_0,x) = \bar p_0\,x_0 + U(x)$, namely the function one differentiates in the extended space, is not proper.}.
  Let us also anticipate that  our result is valid under general hypotheses on $f$ and $l$ that do not guarantee neither uniqueness nor bounded length in finite time of the trajectories (see next subsection).

Investigations on  this kind of value functions have been pursued in several papers, mainly from the point of view of the corresponding Hamilton-Jacobi equation. Indeed  the existence of pairs $(x,a)$ such that $l(x,a)=0$ raises various non trivial problems about  uniqueness. A likely incomplete bibliography,  also containing applications (for instance,   the F\"uller and  shape-from-shading problems),  includes  \cite{BCD}, \cite{I}, \cite{IR}, \cite {CSic}, \cite{Sor2},  \cite{M},   \cite{Ma},  and the references therein.

 Actually, besides displaying an obvious control theoretical meaning (see also Remark  \ref{FC}), our result provides a sufficient condition for the value function to be continuous on  the target's boundary. This continuity property is crucial to recover uniqueness of the corresponding boundary value problem (see Remark \ref{Vcont}).

%\footnote{Because of  condition (\ref{sign}), $  (z,\alpha)\to {\mathcal J}_{z,\alpha}$ might be called a {\it pseudo-time} functional The time-functional is recovered by setting $l\equiv 1$.}.

We conclude this informal presentation  by observing that if some $\eta>0$  existed such that $l(x,a)\geq \eta$ \, $\forall (x,a)\in \R^n\times A$, then
the problem could be easily reduced to an actual  optimal time problem by just utilizing the  reparameterized dynamics
$$\frac{d y}{d\tau} = \ds\frac{f(y,\alpha)}{l(y,\alpha)} .$$
Indeed, after the  (bi-Lipschitz) time-parameter change $\tau(t)\doteq \int_0^t l\,dt$ the Lagrangian turns out to be  transformed into the constant value $1$.
Yet, if only the  weaker sign condition (\ref{sign}) is assumed, a direct approach based on such a  reparameterization cannot be adopted.

\vv

 The paper is organized as follows. In the remaining part of the present  section  we state  rigorously the main result of the paper (Theorem \ref{Th3.1gen}) and provide  some basic definitions. In Section \ref{Heur} we sketch the main result's proof   by heuristic  arguments  and give a geometrical description of the thesis.  Section \ref{Heur} ends with some   examples.  Section \ref{proof} is the longest one, and is entirely devoted to the proof of Theorem \ref{Th3.1gen}, while some technical results are proved in Section \ref{pprop}.   In Section \ref{Concl} we make some  remarks on the meaning of assumption (\ref{C11}) as a viscosity supersolution condition.

\subsection{Precise statement of the main result}
Our main technical assumptions are:

\begin{itemize}
\item[\it (i)] for given positive integers $n$, $m$, the controls $\alpha(\cdot)$ take values in a compact set $A\subset\R^m$ and are Borel-measurable, while the state values $z(t)$, $x$ range over $\R^n$;

\item[\it (ii)] the  {\it target}   ${\bf C}$ is  closed and has   compact boundary;

\item[\it (iii)] the augmented vector field $(l,f)$ is a continuous function on $\C^c\times A$.

     \end{itemize}
     In particular, for any given control $\alpha(\cdot)$ and initial condition $z(0)=x$, the Cauchy problem associated to  the differential equation (\ref{ode}) may have multiple Carath\`eodory solutions. Moreover, the latter may have unbounded velocity near $\C$, since  $f$ itself  is not assumed to be neither  Lipschitz nor bounded near $\C$.
Actually it may well happen that approaching  trajectories  fail  to  {\it reach} the target even when   ${T}_{z} <+\infty$ (see Example \ref{Ex2}).

%$$H(x,p_0,p)\doteq \inf_{a\in A}\left\langle (p_0, p) \,,\,\Big(l(x,a),f(x,a)\Big)\right\rangle.
%$$
Let us introduce the notion of  {\it Minimum Restraint Function}.

 \begin{Definition}\label{CMRL}   We say that a continuous function $U:\overline{\C^c} \to\R$ is a {\em Minimum Restraint Function,}  in short, a  (MRF),  if  $U$   is locally  semiconcave, positive definite, and
 proper on  ${\C}^c$,  and, moreover,
  there exists $\bar p_0\geq 0$ such that
\begin{equation}\label{C1}
H(x,\bar p_0,D^* U(x))<0
\end{equation}
   holds true for all $ x\in \C^c$. \footnote{We refer to Subsection \ref{sub-i2} for the definitions
 of {\it limiting gradient,}
 $D^*U$, and of {\it proper, positive definite and semiconcave function}.}
\end{Definition}

\vskip0.7truecm

\noindent
\begin{Theorem}\label{Th3.1gen}
  Let  a Minimum Restraint Function $U$ exist. Then:
 \begin{itemize}
 \item[{\bf (i)}] the system (\ref{ode}) is   globally asymptotically controllable\footnote{See Definition \ref{(GAC)}.} to $\C$;
  \item[{\bf (ii)}] if   $U$  is a Minimum Restraint Function with $\bar p_0>0$,  then
  \bel{Vprop}
 {{\mathcal V}}(x)\le \frac{U(x)}{{{\bar p_0}}}\,  \qquad\forall x\in \C^c.
\eeq
\end{itemize}
\end{Theorem}

\vskip 0.3 truecm
\begin{Remark}{\rm Petrov-like inequalities (see Example \ref{exP}) are included in  condition (\ref{C1}). However, let us point out that  they only concern  the case where $l=1$ and the dynamics is bounded near the target,  which implies that  optimal (or quasi-optimal) trajectories take a finite time to reach the target.  Instead when the Lagrangian is just non-negative, condition (\ref{C1}) does not force optimal trajectories to approach the target in finite time.   }\end{Remark}
%Let us notice preliminarily  that even if the  dynamics $f$  is globally bounded it may well  happen that ${\mathcal J}_{z,\alpha}<\infty$ and ${T}_{z}=+\infty$: roughly speaking, in order that ${\mathcal J}_{z,\alpha}$ to be finite it is the magnitude of the {\it ratio} $f/l$ that plays a crucial role. On the other hand,  this ratio could be neither everywhere defined (since $l$ is allowed to be zero) nor bounded (even in the case $l=1$)

\vskip0.3truecm

\begin{Remark} {\rm
Because of the sign assumption (\ref{sign}),  a Minimum Restraint Function, (MRF),  is in particular a Control Lyapunov Function, (CLF).  Hence,  as a byproduct of statement {\bf (i)} in Theorem \ref{Th3.1gen} we get an extension to merely continuous, unbounded dynamics of the results concerning the relation between  Control Lyapunov Functions  and asymptotic controllability (see e.g. \cite{S2}).
}
\end{Remark}

\begin{Remark}\label{Vcont} {\rm
 Because of the bound  (\ref{Vprop}),  the  value function ${{\mathcal V}}$
%inherits some regularity of $U$  on the boundary ${\partial \bf C}$ of the target  ${\bf C}$. ${{\mathcal V}}$
turns out to be continuous on ${\partial \bf C}$. Actually, this is a theoretical motivation for a result like
Theorem \ref{Th3.1gen}, in that  the continuity of the value function on the target's boundary is essential to establish   comparison, uniqueness,
%\footnote{ The drawback caused by the possibility that $l=0$ at some point, has a classical PDE's counterpart in the fact that some points $(x,p_0,p)$ might be {\it characteristic} for the Hamiltonian
%$\hat H(x_0,x,p_0,p)\doteq H(x,p_0,p)$
%and the hyperplanes $\{(x_0,x)\in\R^{1+n}, \,\,tx_0=constant \}$, which in our case means $\ds\frac{\partial \hat H}{\partial p_0} = 0$. As suggested by the theory of characteristics,  this may cause  non-uniqueness of solution for a Cauchy problem  associated to the Hamilton-Jacobi equation $\hat H(x,D_{(x_0,x)}u) = 0$.  },
and robustness properties  for the associated Hamilton-Jacobi-Bellman equation (see \cite{MS}, \cite{Sor1}, and \cite{M}).}
\end{Remark}

\begin{Remark}\label{FC}{\rm  Our investigation might be useful also for feedback control. Indeed, on the one hand  the uniqueness of the H-J equation is an obvious ingredient in a feedback-oriented construction. On the other hand, a (MRF) $U$ can be possibly exploited in order  to build a "safe" feedback law for problem (\ref{P}) (see e.g. \cite{CLSS} and  \cite{AB} for general reference to  feedback stabilization).
}\end{Remark}

\begin{Remark} {\rm
It is   easy to adapt Theorem \ref{Th3.1gen} to the case when the state space is an open set $\Omega\subset\R^n$, $\Omega\supset\C$.  In fact,  the thesis keeps unchanged as soon as  one requires the (MRF) $U:\Omega\setminus\overset{\circ} \C\to\R$ to verify all the assumptions in  Definition \ref{CMRL} in $\Omega$,   plus the following one:
$$
\exists U_0\in]0,+\infty]: \ \
\lim_{x\to x_0, \  x\in\Omega}U(x)=U_0 \ \  \forall x_0\in\partial\Omega;  \quad
U(x)<U_0 \quad\forall x\in\Omega\setminus\overset{\circ} \C.
$$}
\end{Remark}

\subsection{Basic definitions}\label{sub-i2}

For the reader convenience, some classical concepts, like  (GAC), and
a few technical definitions (part of which have already been used in Theorem \ref{Th3.1gen} above) are here recalled.

\begin{Definition} {\rm (Positive definiteness).} A continuous function $U:\overline{\C^c} \to\R$ is said {\rm positive definite on $\C^c$} if  $U(x)>0$ \,$\forall x\in{\C}^c$ and $U(x)=0$ \,$\forall x\in\partial{\bf C}$. Moreover $U$ is called {\rm proper  on ${\C}^c$} if $U^{-1}(K)$ is compact as soon as  $K\subset[0,+\infty[$ is compact.
\end{Definition}

 \begin{Definition} {\rm (Semiconcavity). }\label{sconc} Let $\Omega\subset\R^n$ be an open set, and let  $F:\Omega\to\R$ be a continuous function.  $F$is said to be {\rm locally semiconcave  on $\Omega$} if for any point $x\in\Omega$ there exist $R>0$ and $\rho>0$ such that
 $$
 F(z_1)+F(z_2)-2F\left(\frac{z_1+z_2}{2}\right)\le R|z_1-z_2|^2 \quad \forall z_1, \, z_2\in B_n(x,\rho).
 $$
 \end{Definition}
Let us remind that locally semiconcave functions are locally Lipschitz. Actually, they are twice differentiable almost everywhere (see e.g. \cite{CS}).
 \begin{Definition}\label{D*}{\rm (Limiting gradient). }  Let $\Omega\subset\R^n$ be an open set, and let  $F:\Omega\to\R$  be a locally Lipschitz function.  For every $x\in \Omega$ let us set
$$
D^*{F}(x) \doteq \Big\{ w\in\R^n: \ \  w=\lim_{k}\nabla {F}(x_k), \ \  x_k\in DIFF(F)\setminus\{x\}, \ \ \lim_k x_k=x\Big\}
$$
where   $\nabla$ denotes the classical gradient operator and $DIFF(F)$ is the set of differentiability points of $F$.  $D^*{F}(x)$ is called
the  {\rm set of limiting gradients} of $F$ at $x$.
\end{Definition}
For every $x\in\Omega$,   $D^*{F}(x)$ is a  nonempty, compact subset of $\R^n$ (more precisely, of the cotangent space $T^*_x\Omega$). Notice that, in general,  $D^*{F}(x)$ is not convex\footnote{Actually its  convexification coincides with the Clarke's generalized gradient.}.

\vskip0,4truecm
To give the notion of global asymptotic  controllability we need to recall the concept of {\it function belonging to ${\cal KL}$}:  these are continuous functions  $\beta:[0,+\infty[\times[0,+\infty[\to[0,+\infty[$ such that: (1)\, $\beta(0,t)=0$ and $\beta(\cdot,t)$ is strictly increasing \footnote{We call  a real map $\phi$  {\it decreasing [increasing],}   if $\phi(r_1)\geq \phi(r_2)$ [$\phi(r_1)\leq \phi(r_2)$ ] as soon as $r_1<r_2$ and  {\em strictly  decreasing [increasing],}  if the inequality is always strict.}  and unbounded for each $t\ge0$; (2)\, $\beta(r,\cdot)$ is decreasing  for each $r\ge0$; (3)\, $\beta(r,t)\to0$ as $t\to+\infty$ for each $r\ge0$. For brevity, let us use the notation ${\bf d}(x)$ in place of ${\bf d}(x,\C)$.

 \begin{Definition}\label{(GAC)}\footnotemark{\rm(GAC)} The system   (\ref{ode}) is   {\rm globally asymptotically controllable to $\C$}  --shortly, (\ref{ode}) is {\rm (GAC) to $\C$}--  provided  there is a  function $\beta\in{\cal KL}$  such that,  for each initial state $x\in\C^c$,    there exists an admissible trajectory-control pair $(z,\alpha):[0,+\infty[\to\R^n\times A$   that verifies
\bel{bbound}
{\bf d}(z(t))\le \beta({\bf d}(x),t) \qquad \forall t\in[0,+\infty[.
\eeq
\end{Definition}
 \footnotetext{
 To be precise, we are considering a slight variation  of the standard notion of {\rm(GAC)} to $\C$, which would require  $\C$ to be weakly invariant with respect to the control dynamics, since we are interested in the behavior of any admissible trajectory $z$ just for $t\in[0,T_z[$.  Therefore, we  fix an arbitrary $\bar z\in\partial\C$ and,  when $T_z<+\infty$,  we prolong  $z$ to $[0,+\infty[$ by setting $z(t)=\bar z$ for all $t\ge T_z$. }
 \vskip0.3truecm
 In  Theorem \ref{Th3.1gen}'s  proof we shall make  use of the notions of  {\it partition} of an interval and of its {\it diameter}. To avoid vagueness let us state the precise meaning we attach to these terms.

 \begin{Definition}\label{pi}   Let us consider an interval $[0,b[$, \, $b\in]0,+\infty]$.  A  {\it partition}  of $[0,b[$ is a sequence $\pi=(t^j) $ such that
$t^0=0, \quad t^{j-1}<t^j$ \, $\forall j\ge 1$, and
  {\rm either}  $\lim_{j\to+\infty}t^j=b$ {\rm or}  there exists some $n_\pi\in\N$ such that $t^{n_\pi}=b$. In the latter case, we say that $\pi$ is a {\rm finite partition} of $[0,b[$.
  The number diam$(\pi)\doteq\sup_{\{j: \ t^j\le b\}}(t^{j }-t^{j-1})$ is called the {\rm diameter}  of the sequence $\pi$.
  \end{Definition}

\section{Heuristics  of the proof and some examples}\label{Heur}

\subsection{A one dimensional differential inequality}
As in the case  of asymptotic controllability, the underlying idea of Theorem \ref{Th3.1gen} relies on a one-dimensional argument involving a differential inequality. To express this issue, let us assume  some simplifying facts. Let us  begin with making the hypothesis that $U$ is of class $C^1$, so that   assumption (\ref{C1})  reads
\bel{Csimple}
\inf_{a\in A}\left\langle (\bar p_0, \nabla U(x)) \,,\,\Big(l(x,a),f(x,a)\Big)\right\rangle<0 \qquad \forall x\in{\bf C}^c.
\eeq
Let us also assume that there exists a continuous selection $$a(x)\in A(x)\doteq\left\{a\in A: \  \left\langle (\bar p_0, \nabla U(x)) \,,\,\Big(l(x,a),f(x,a)\Big)\right\rangle<0\right\},$$
 so that (\ref{Csimple}) yields
\bel{senzainf}
\left\langle (\bar p_0, \nabla U(x)) \,,\,\Big(l(x,a(x)),f(x,a(x))\Big)\right\rangle<0  \qquad \forall x\in{\bf C}^c
\eeq
for some $\bar p_0\ge0$.

Let us consider a  solution of the augmented  Cauchy problem
\bel{feedback}
\left\{\begin{array}{l} \dot z_0 = l(z,a(z))\\\,\\\dot z = f(z,a(z))\\\,\\
(z_0,z)(0) = (0,x)\end{array}\right.
\eeq
and let us set
$$\xi(t)\doteq U(z(t)).$$
Then, by (\ref{senzainf}),
\bel{odi1}
\begin{array}{l}
\dot \xi(t) + \bar p_0\, \dot z_0(t)= \langle \nabla U(z(t)) , \dot z(t) \rangle + \bar p_0\, \dot z_0(t) =  \\ \, \\
\qquad\qquad\langle \nabla U(z(t)) \,,\, f(z(t), a(z(t))\rangle + \bar p_0\, \dot z_0(t)  < 0
\end{array}
\eeq
Notice that $\xi(0) = U(x) >0$ and $z_0(0)=0$. The rough idea of the proof (to be sharpened  through suitable nonsmoothness' and o.d.e.'s
arguments) amounts to show that:\begin{itemize}
\item[(A)] $\xi(t)$ is defined, strictly decreasing and  tends to zero in a possibly unbounded interval $[0, T[$ : this means that $\lim_{t\to T}{\bf d}(z(t),\C)=0$, which coincides with part {\bf (i)} of the thesis of Theorem \ref{Th3.1gen}.
    \item[(B)] If $\bar p_0>0$, the rate of growth of the (non negative and increasing) map $z_0(t)$ is bounded by  $\ds -\frac{\dot\xi(t)}{\bar p_0}$: this implies that
        $$\mathcal V(x)\leq \lim_{t\to T} z_0(t) \leq  \frac{1}{\bar p_0}\left( \xi(0) -  \lim_{t\to T} \xi(t)\right) = \frac{U(x)}{\bar p_0},$$
       which coincides with the statement {\bf (ii)} of Theorem  \ref{Th3.1gen}.
\end{itemize}

Let us observe, however, that inequality (\ref{odi1}) is not enough to prove (A). In fact, $\xi(\cdot)$ could well decrease asymptotically to a value $\bar \xi>0$. To show that  $\xi(\cdot)$ actually converges  to zero, one  deduces  the differential inequality
\bel{odi2}
\dot \xi(t) + \bar p_0\, \dot z_0(t) = \leq - m(\xi(t)),
\eeq
from (\ref{odi1}), $m(\cdot)$ being a suitable  positive,
strictly increasing function on $]0,+\infty[$.  This is, in fact, the essential content of Proposition \ref{claim1}, where (\ref{odi2}) is replaced by the nonsmooth relation (\ref{Hg}). The other ingredient of the proof is Proposition \ref{claim2bis}, where nonsmooth analysis techniques are applied to show that
  things actually work even  without the simplifying regularity we are assuming here. In particular, when $\bar p_0>0$ one can emulate the above (B) and get the bound on the value function.

\subsection{A geometrical insight}
A further interpretation of the result in part (ii) of Theorem \ref{Th3.1gen}  in the case $\bar p_0>0$,
is provided by the following "geometrical" description of the above heuristic  arguments. To begin with,  notice that the {$(1+n)$-dimensional target $\tC= [0,+\infty[\times\C$ has no longer compact boundary.  Therefore, no proper (CLF) can exist. Actually, a (MRF) $U$,  when considered as a function on   $[0,+\infty[\times \R^n$, is {\it not proper}. On the other hand,   let us consider the map
$$
W(x_0,x)\doteq U(x) +\bar p_0\, x_0.
$$

The key point which makes Theorem \ref{Th3.1gen}  work relies on the following  three facts:
\begin{itemize}\item[1)] $W$ is proper in $[0,+\infty[\times \R^n$;
 \item[2)] the inequality
\bel{liapk}
H(x, D^*_{(x_0,x)} W(x_0,x))<0,
\eeq
which coincides with inequality (\ref{C1}), says that, for every $ \bar x\in \C^c$,  there exist trajectories of the augmented system $(\dot z_0,\dot z)=(l,f)$  starting from $(0,\bar x)$ and  remaining  inside {\it the $\bar x$-influence set}
$$
 \Big\{ (x_0,x)\in[0,+\infty[\times \R^n: \quad W(x_0,x)\leq U(\bar x)\Big\};
$$
\item[3)] the level sets of $W$ intersect the extended target $\tC$.
\end{itemize}
The situation is illustrated in Fig.1.
\newpage
\begin{figure}
\centering
\includegraphics[scale=0.30]{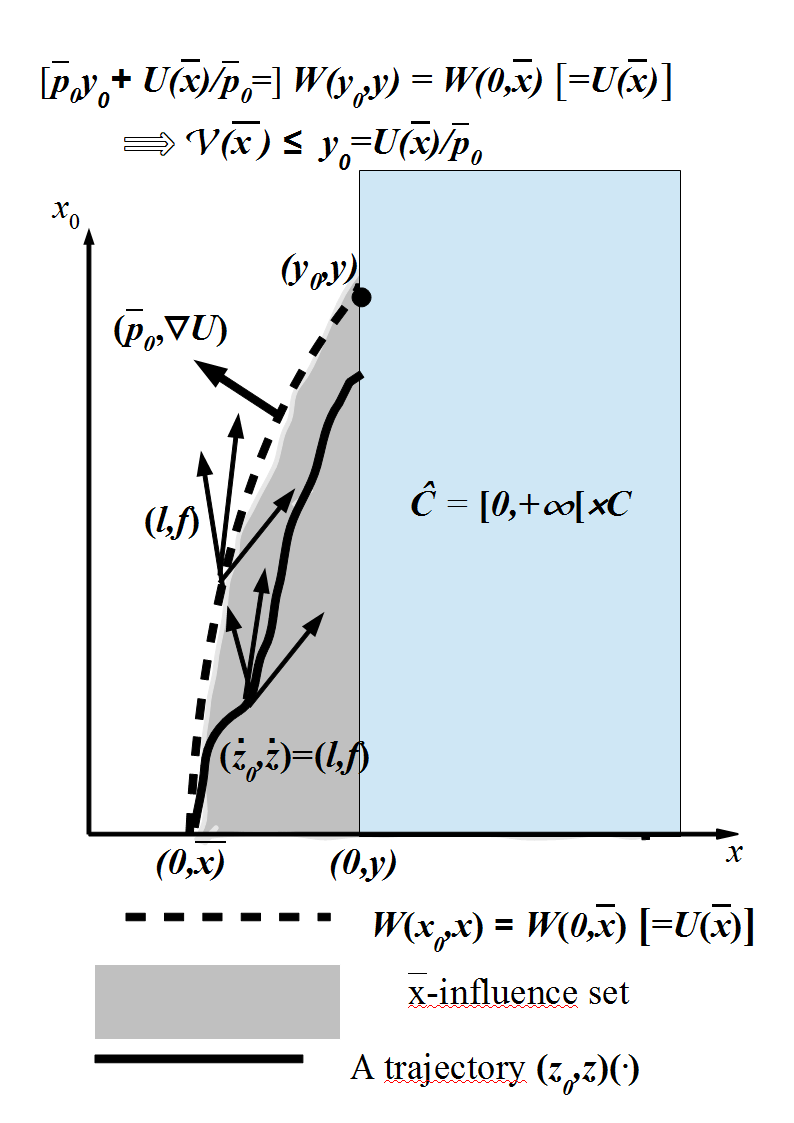}
\caption{The   level sets of  $W$  intersect the extended target $\tC$}
\end{figure}

\subsection{Some examples} \label{Ex}
\begin{Example}\label{exP}{\rm A prototype of (MRF) is a map of the form $ U(x)=\Phi\circ{\bf d}(x)$, where
 $\Phi:[0,+\infty[\to\R$ is a continuous map such that   $\Phi(0)=0$ and its restriction to  $]0,+\infty[$ is a strictly increasing $C^1$-diffeomorphism.

\noindent  In particular,  in the minimum time problem, where $l\equiv 1$,  the  inequality   (\ref{C1})  includes   the following  {\it weak Petrov condition} (see e.g. \cite{S1},  \cite{CS},  and \cite{BP}):
\begin{itemize}
\item[ {\bf(P)}] there exist $\delta>0$ and a continuous, increasing map $\mu:[0,\delta]\to[0,+\infty[$  verifying $\mu(0)=0$, $\mu(\rho)>0$ for $\rho>0$,  $\ds\int_0^\delta\frac{d\rho}{\mu(\rho)}<+\infty$, and such that,  $\forall x\in\C^c$ with ${\bf d}(x)<\delta$,  one has
\bel{genPet}
\inf_{a\in A}\langle D^*{\bf d}(x), f(x,a)\rangle\le -\mu({\bf d}(x)).
\eeq
\end{itemize}
Indeed, if we set   $\Phi(r)\doteq  \int_0^r\frac{d\rho}{\mu(\rho)}$ for all $r\in[0,\delta]$ and choose an arbitrary $\bar p_0\in[0,1[$, we can write (\ref{genPet}) as
\footnote{This equivalence follows from the straightforward  set identity $D^*(\Phi\circ{\bf d})(x)=\nabla \Phi({\bf d}(x))\, D^*{\bf d}(x)$}
$$
H(x,{\bar p_0},D^*(\Phi\circ{\bf d})(x))=\inf_{a\in A}\left\{{\bar p_0} + \left\langle D^*(\Phi\circ{\bf d})(x), f(x,a)\right\rangle\right\}\le -(1-{\bar p_0})<0.
$$
In particular,  whenever $\Phi$ is   linear one recovers the classical Petrov condition.
}
 \end{Example}

\begin{Example}\label{Exa1}{\rm
 Let $r,s$ be arbitrary real numbers and  let $\psi_1(x)$, $\psi_2(x)$ be continuous  functions such that $\psi_1(x)\ge M_1$ and $0\le \psi_2(x)\le M_2$ for some   $M_1$, $M_2>0$.   Consider  the target $\C\doteq\{0\}$, the  control dynamics
 $$
\dot z = a\,\psi_1(z) |z|^r \qquad a\in\{-1,1\}
 $$
 and the  Lagrangian
$$
l(z)\doteq \psi_2(z)|z|^{s}.
$$
Notice that   $l$ can be zero on an arbitrary subset of its domain.

Let  $s,r$ verify
\bel{F1}
 s-r>-1
\eeq
 and consider the map
$$
 U(x)\doteq \frac{M_2}{M_1(s-r +1)}\,|x|^{s-r +1}.
 $$
Observe that $U$ is   proper, positive definite and semiconcave  on $\R\setminus\{0\}$ ---actually, $U\in C^\infty(\R\setminus\{0\})$.  Moreover, for any $p_0\ge0$,
$$
\begin{array}{c}
H(x, p_0,p)  = \inf_{a\in\{ -1, 1\}}\Big\{ p_0  \psi_2(x)\,|x|^{s} + ap\, \psi_1(x)\, |x|^r \Big\} \le\\\,\\
p_0M_2\,  |x|^{s}-M_1\,|p||x|^r.
\end{array}
$$
Since
$$
D^*U(x) = \left\{\nabla U(x)\right\} = \left\{ \text{sign}(x)\,\frac{M_2}{M_1\,} |x|^{s-r}\right\} \qquad \forall x\ne 0,
$$
 one obtains
$$
H(x,{\bar p_0},D^*U(x)) <0 \qquad \forall x\neq 0,
$$
as soon as ${\bar p_0}\in ]0,1[$.

By applying Theorem \ref{Th3.1gen}  we get that the  value function ${\mathcal V}$ ---namely, the minimum cost to attain the target (possibly in infinite time) ---, verifies
$$
{\mathcal V}(x) \leq \frac{M_2}{M_1 \bar p_0(s-r +1)}\,|x|^{s-r +1}\qquad\forall x\in \R
$$
for all  $\bar p_0\in]0,1[$, which implies
$$
{\mathcal V}(x) \leq \frac{M_2}{M_1  (s-r +1) }\,|x|^{s-r +1} \qquad\forall x\in \R.
$$

Notice that (\ref{F1}) is crucial. Indeed, if $s-r\leq -1$,  a nonsingular (MRF) { \it may fail to exists}. In  fact, consider the trivial case when $\psi_1=\psi_2\equiv  1$. If $U$ were a nonsingular (MRF),   for almost every $x\in]0,1]$ we should have
$$
|\nabla U(x)| >\bar p_0\, x^{s-r}\geq \bar p_0\, x^{-1}
$$
for some $\bar p_0>0$, which clearly  prevents $U$ to be positive definite.}
\end{Example}

\if{
\begin{Remark}\label{Exa 2} {\rm  The function  $U$ in  Example \ref{Exa1} is  of class $C^\infty$ (in the interior of its domain). This is more the  exception rather than the rule. Actually,  it is sufficient to slightly modify the previous example  to get a situation where a (MRF) cannot be   smooth. This would be the case, for instance, when $\C=\{0,1\}$, $
f(z,a) = a\,  \psi_1(z){\bf d}^r(z),\quad a\in \{-1,1\}
$,
and $
l(z)\doteq \psi_2(z) {\bf d}^s(z).
$
Indeed, if $\tilde U$ were a smooth (MRF), it should attain a maximum at some point $\tilde x \in]0,1[$, so that one would have $D^*\tilde U(\tilde x) = \{\nabla\tilde U(\tilde x)\} =  \{0\}$. Since $H(\tilde x,p_0,0)  =p_0\,l(\tilde x) \geq 0$, for any $p_0\geq 0$  condition (\ref{C1}) would fail to  be verified at $\tilde x$.
}

  \end{Remark}
  \fi
\vskip0.4truecm
In Example \ref{Exa1}  the dynamics may happen to be unbounded (precisely, when $r<0$). However, when the time $T_z$ to approach the target happens to be finite, the trajectory's interval can be prolonged to the closed interval $[0,T_z]$, so that  $z(T_z)\in\C$.  Let us remark that this is due to the one-dimensionality of the state space. In fact in the next example we see, among other things, that there is a connected component of the target that can be approached in finite time while {\it it cannot be reached.} \footnote{Incidentally, this is the reason why we have adopted a notion of (GAC) slightly more general than the usual one.}
 
\begin{Example}\label{Ex2}{\rm Let us set
$$\C\doteq \C_1\cup\C_2\subset \R^2,$$
where
$$\C_1\doteq \{x\in\R^2, \,\,\,|x|\leq 1\},\qquad \C_2\doteq \{x\in\R^2, \,\,\,|x|\geq 4\}, $$
and let us consider the control dynamics
$$\dot z=\frac{M\cdot z}{(|z|-1)} - \alpha z$$
with the control $\alpha$ taking values in $\{-1,1\}$ and $M\doteq \begin{pmatrix} \,\,0  & 1\\-1& 0\end{pmatrix}
$.
Finally let us consider the Lagrangian
$$
l(x)  \doteq \left\{\begin{array}{ll} (|z|-1)^2k(z)\qquad &\hbox{if} \,\,1\leq |z|\leq 2\\(3-|z|)^2k(z)\qquad &\hbox{if} \,\,2<|z|\leq 3\\
0  \qquad &\hbox{if}\,\, 3<|x|\leq 4\end{array}\right.,
$$
where $k:\C^c\to\R$ is any continuous function verifying $0\leq k(z)\leq 1\,\,\,\,\forall z\in \C^c $.
Let us  begin with showing that the target's component $\C_1$ can be indefinitely approached in finite time but {\it it cannot be reached}. Indeed, in polar coordinates $z=\rho e^{i\theta}$ the control equation becomes
$$
\left\{\begin{array}{l}
\dot \rho = - \alpha\rho\\\,\\
\dot\theta = -(\rho -1)^{-1}\end{array}\right.
$$
 Let $\bar\rho\in]1,4[$, $\bar\theta\in[0,2\pi[$.  Choose $\bar z=\bar\rho\, e^{i \bar\theta}\in\C^c$,  $\alpha\equiv 1$, and consider  the solution of  the first equation: $\rho:[0,\ln\bar\rho[\to\R$, $\rho(t)\doteq \bar\rho\, e^{-t}$. The time this trajectory takes to approach $\C$ is equal to $\ln \bar\rho$. Furthermore,
  $$
\lim_{t\to \ln \bar\rho} (\rho(t),\theta(t)) =(1,+\infty).
$$
It follows that
$$
\lim_{t\to \ln\bar\rho}{\bf d}( z(t)) =0
$$
{\it but $z(t)$ has no limit} as $t\to \ln\bar\rho$ . Namely, the trajectory spirals faster and faster around $\C_1$ while approaching it.

On the other hand it is trivial to see that the target's component $\C_2$ {\it can be reached} in finite time by implementing the constant control $\alpha\equiv-1$.
Lastly, setting for every $\epsilon\geq 0$,
$$
U_\epsilon(z)  \doteq \left\{\begin{array}{ll} 2\epsilon  + \frac{(|z|-1)^3}{3}\qquad &\hbox{if} \,\,1\leq |z|\leq 2\\\,\\\ds \epsilon(4-|z|) + \frac{(3-|z|)^3}{3}\qquad &\hbox{if} \,\,2<|z|\leq 3\\\,\\\ds
\epsilon(4-|z|)\qquad &\hbox{if} \,\,3<|z|\leq 4\end{array}\right.,
$$
 one can easily check that as soon as $\epsilon>0$ $U_\epsilon$ is a  Minimum Restraint Function with $\bar p_0 = 1$. Therefore, if $\V(x)$  denotes the value function of the problem, by Theorem \ref{Th3.1gen} we get
  $
  \V(z)\leq U_\epsilon(z)$ for every $\epsilon>0$, so that
 \bel{MRFex}
  \V(z)\leq U_0(z).\eeq

  Notice  that, in agreement with the definition of  (MRF),  the functions $U_\epsilon$ are not smooth at their maximum points. Moreover, the inequality   (\ref{MRFex}) is optimal  in the ring  $R\doteq\{ z\,\,|\,\,3\leq |z| \leq 4\}$, in that $U_0=0$, and hence $\V=0$, on $R$.

}
\end{Example}\vskip1truecm

\section{Proof of the main result}\label{proof}
\subsection{Preliminary results}

The proof of Theorem \ref{Th3.1gen} relies on   Propositions \ref{claim1} and  \ref{claim2bis} below. In order to retain clarity in  the main proof's  argument, we postpone the proofs of these technical results to Section \ref{pprop}.

\begin{Proposition}\label{claim1} Let $U$ be a  (MRF)  and let ${\bar p_0}\ge0$ make  (\ref{C1}) hold true, i.e.
$$
H(x,{{\bar p_0}},D^* U(x))<0 \qquad \forall x\in{\C}^c.
$$
  Then for every $\sigma>0$ the map $U$ verifies also the differential inequality
\bel{Hg}
 H(x,{{\bar p_0}},D^* U(x))\le -m(U(x)) \qquad \forall x\in U^{-1}(]0,\sigma[),
\eeq
where  $m:[0,+\infty[\to\R$ is a suitable  continuous, strictly increasing  function  verifying $m(r)>0$ \, $\forall r>0$.
\end{Proposition}

\begin{Remark}\label{glob} {\rm
If we replace condition (\ref{C1}), namely
$$
\exists {\bar p_0}: \ \ H(x,{{\bar p_0}},p)<0 \qquad \forall x\in\C^c, \ \ \forall p\in D^* U(x),
$$
  with the stronger assumption $$
\exists {\bar p_0}: \ \  \forall M>0,\quad \sup_{\, {\bf d}(x)\ge M; \ p\in  D^*U(x)} H(x,{\bar p_0},p)<0,
$$
 then   by slight changes in  the proof of Proposition \ref{claim1} below, one can prove that there exists a continuous, strictly increasing  function  $m :[0, +\infty[\to\R$, independent of $\sigma$,  verifying $m(r)>0$ \, $\forall r>0$, such that (\ref{Hg}) holds  for all $x\in\C^c$.}
\end{Remark}

\begin{Proposition}\label{claim2bis}
Let $U$ be a  (MRF) and let $\sigma>0$.   Let  $m$ be  defined as in Proposition \ref{claim1}. Fix $\varepsilon>0$ and   $\bar\mu$,  $\hat\mu$    such that   $0<\hat\mu<\bar\mu<\sigma$. Then there is some $\delta>0$ such that, for every $\delta'\in]0,\delta[$ and  for  each $x\in \C^c$ verifying $U(x)=\bar\mu$,  for a suitable $\hat t>0$ one can construct  a trajectory-control pair
$$
(z,\alpha):[0,\hat t]\to U^{-1}([\hat\mu,\bar\mu])\times A, \qquad z(0)=x,
$$
and a finite partition $\pi=\{t^0,\dots,t^{\bar n}\}$ of $[0,\hat t]$ such that diam$(\pi)\le\delta'$ and  with  the following properties:
\begin{itemize}
\item[\bf (a)]  $U(z(\hat t))=\hat\mu<U(z(t))\le U(z(0))=U(x)=\bar\mu$ \, $\forall t\in[0,\hat t[$;
\item[\bf (b)]  for  every $j\in\{1,\dots,\bar n\}$ and $\forall t\in[t^{j-1},t^j[$,
\bel{stb}
U(z(t^j))<U(z(t))\le U(z(t^{j-1}));
\eeq
\item[\bf (c)]  for  every $j\in\{1,\dots,\bar n\}$ and $ \forall  t\in[t^{j-1}, t^j]$,
\bel{dainserire}
\begin{array}{l}
U(z(t))-U(z(t^{j-1}))+ \frac{\bar p_0}{\varepsilon+1}\int_{t^{j-1}}^t l(z(t),\alpha(t))\, dt \le \\ \, \\
\qquad\qquad\qquad\qquad\qquad\qquad  -\frac{1}{\varepsilon+1}\int_{t^{j-1}}^t m(U(z(t)))\, dt.
  \end{array}
\eeq
\end{itemize}
\end{Proposition}

\vv

\subsection{Proof of Theorem   \ref{Th3.1gen}.}

In the sequel, given a constant $\mu>0$,  for any continuous path $y:[\tau, +\infty[\to\R^n$ with  $U(y(\tau))>\mu$,  we   define the time to reach the enlarged target $U^{-1}([0,\mu])$ as
\bel{Tz}
{\cal T}_y^\mu\doteq\inf\{r\ge\tau: \ U(y(r))\le \mu\}
\eeq
(${\cal T}_y^\mu=+\infty$ if $U(y(r))> \mu$ for all $r\ge\tau$).

\vv
 Let   $\sigma>0$ be a positive constant and  let $m$ be defined as in Proposition \ref{claim1}. Fix $\varepsilon>0$ and let $(\nu_k)\subset]0,1]$ be   a sequence   such that $1=\nu_0>\nu_1>\nu_2>\dots$ and $\lim_{k\to\infty}\nu_k=0$.   Assume that $x\in U^{-1}(]0,\sigma[)$ and set
$$
\mu_k\doteq \nu_k U(x) \quad \forall k\ge 0.
$$

We are  going to exploit   Proposition \ref{claim2bis} in order to build  the trajectory-control pair
$$
(z,\alpha):[0,\bar t[\to\C^c\times A
$$
by   concatenation,
$$
(z(t),\alpha(t)) = (z_k(t),\alpha_k(t)) \quad  \forall t\in [t_{k-1},t_k[, \quad \forall k\ge1.
$$

{\it Step $k=1$.} Let us begin by constructing  $(z_1,\alpha_1)$.  Setting    $\bar\mu=\mu_0$,  $\hat\mu=\mu_1$,     let $(z_1,\alpha_1):[0,\hat t]\to U^{-1}([\mu_1,\mu_0])\times A$ be a trajectory built according to  Proposition \ref{claim2bis} such that $z_1(0)=x$.  We set $t_0\doteq 0$ and $t_1\doteq \hat t$ and we observe that $t_1= T_{z_1}^{\mu_1}$, in view of {\bf (a)} in Proposition \ref{claim2bis}.

{\it Step $k>1$.}  Let us proceed by defining $(z_k,\alpha_k)$ for $k> 1$. Setting  $\bar\mu=\mu_{k-1 }$,  $\hat\mu=\mu_{k }$,     let $(\hat z_k,\hat \alpha_k):[0,\hat t]\to U^{-1}([\mu_k,\mu_{k-1}])\times A$ be a trajectory built according to  Proposition \ref{claim2bis} such that $\hat z_k(0)=z_{k-1}(t_{k-1})$. We set  $t_k\doteq t_{k-1}+\hat t$ and  $(z_k,\alpha_k)(t)=(\hat z_k,\hat \alpha_k)(t-t_{k-1})$ $\forall t\in[t_{k-1},t_k]$. We observe that $t_k= T_{z_k}^{\mu_k}$, still in view of {\bf (a)} in Proposition \ref{claim2bis}.

The concatenation procedure is concluded as soon as we set   $\bar t \doteq  \lim_{k\to \infty} t_k$. Notice that it may well happen that $\bar t=+\infty$.  We claim that
\bel{raggiunge}
\lim_{t\to\bar t^-}  {\bf d}(z(t)) = 0.
\eeq
 For every $k\ge1$,  let us apply Proposition \ref{claim2bis}, which yields the existence of a finite partition $\pi_k=\{\hat t^0_k,\dots,\hat t^{\bar n_k}_k\}$ of $[0,t_k- t_{k-1}]$ such that, setting,
$$t_k^j\doteq  t_{k-1}+\hat t^j_k  \qquad\forall j\in\{0,\dots, \bar n_k\},
$$
one has      $z(0)\,(=z_1(0))=x$,  and,   for every $k\ge1$:
\begin{itemize}
\item[\bf (a)$_k$]
$z_{k+1}(t_{k}) = z_{k}(t_{k})$, \, $U(z_k(t_{k-1})) = \mu_{k-1 }$;
\item [\bf (b)$_k$] for all $j\in\{1,\dots, \bar n_k\}$, \newline
   $U(z_k(t_k^j))<U(z_k(t))\le U(z_k(t_{k}^{j-1}))\le U(x)$   \,  $\forall t\in[t_{k}^{j-1},t_k^j[$;
\item[\bf (c)$_k$]  for all $j\in\{1,\dots, \bar n_k\}$, \newline
$U(z_k (t))-U(z_k(t_k^{j-1}))+ \frac{\bar p_0}{\varepsilon+1} \int^{t }_{t_k^{j-1}} l(z_k(t),\alpha_k(t))\, dt\le$ \newline
\hphantom{mmmmmmmmmmmmmm} $  -\frac{1}{\varepsilon+1} \int^{t }_{t_k^{j-1}} m(U(z_k(t)))\, dt$   \, $\forall   t\in[t_k^{j-1}, t_k^j]$.
 \end{itemize}

Indeed, in view of point {\bf (b)$_k$} above, (\ref{raggiunge}) is equivalent to
\bel{discrete}
 \lim_{k\to\infty} {\bf d}(z_k(t_{k}))=0.
\eeq
On the other hand,  since  $U$ is proper and  positive definite,  (\ref{discrete}) is a straightforward consequence  of
$$\lim_{k\to\infty} U(z_k(t_k)) = \lim_{k\to\infty} \nu_{k }\,U(x) = 0$$

Therefore (\ref{raggiunge}) is verified.

Notice that  {\bf (b)$_k$}  implies also that
\bel{sigma}
U(z(t))\le U(x)<\sigma \quad \forall t\in[0,\bar t[.
\eeq
\vv

  In order to conclude the proof that the system is (GAC) to $\C$ (part  {\bf (i) } of the theorem),   we have to establish the existence of  a ${\cal KL}$ function $\beta$ as in  Definition \ref{(GAC)}.

   Let $k\ge1$.  From condition {\bf (c)$_k$} and in view of  the definition of $(z_k,\alpha_k)$, we have  $\forall   t\in[t_ {k-1}, t_k^j]$,
$$
\begin{array}{l}
U(z_k(t))-U(z_k(t_{k-1}))
=[U(z_k(t))-U(z_k(t_k^{j-1}))]
+ [U(z_k(t_k^{j-1}))-U(z_k(t_k^{j-2}))]  \\ \, \\  \qquad +\dots+[U(z_k(t_k^{1}))-U(z_k(t_k^{0}))]\le  \\ \, \\
\qquad\qquad\qquad\qquad -\frac{\bar p_0}{\varepsilon+1} \int^{t }_{t_{k-1}} l(z_k(\tau),\alpha_k(\tau))\, d\tau
-\frac{1}{\varepsilon+1} \int^{t }_{t_ {k-1}} m(U(z_k(\tau)))\, d\tau,
 \end{array}
 $$
 which implies that $\forall t\in[0,t_k^j]$,
\bel{EP}
U(z(t))-U(x) \le  -\frac{\bar p_0}{\varepsilon+1} \int^{t }_{0} l(z(\tau),\alpha(\tau))\, d\tau  -\frac{1}{\varepsilon+1} \int^{t }_{0} m(U(z(\tau)))\, d\tau.
\eeq
Being $l\ge0$, in particular  we have
$$
U(z (t ))-U(x)  \le  -\frac{ \int^{t }_{0} m(U(z(\tau))\, d\tau}{\varepsilon+1}.
$$
Since $\min\{m(U(z(\tau)): \ \tau\in[0,t_k^j]\}=m(U(z(t_k^j))$, we get
\bel{stg}
U(z(t ))+\frac{ m(U(z(t_k^j))\,t}{\varepsilon+1}\le U(x) \qquad \forall t\in[0,t_k^j].
\eeq
Observe that the function $\tilde m:[0,+\infty[\to[0,+\infty[$  defined  by $\tilde m(r)\doteq \min\{r,m(r)\}$   for all $r\in[0,+\infty[$  is  continuous,  strictly increasing, and $\tilde m(r)>0$ \, $\forall r>0$,  $\tilde m(0)=0$. Then,  for any $k\ge 1$ and for any $j\in\{1,\dots,{\bar n_k}\}$,
$$
\tilde m(U(z(t_k^j))\left[1+\frac{ t_k^j}{\varepsilon+1}\right]\le U(x),
$$
so that
$$
U(z(t_k^j))\le \tilde m^{-1}\left(\frac{\varepsilon+1}{\varepsilon+1+t_k^j}\,U(x)\right).
$$
Let   $t\in[0,\bar t[$. Then $t\in[t_k^{j-1}, t_k^j[$ for some  $k\ge1$ and some $j\in\{0,\dots,{\bar n_k}\}$.   Moreover, by possibly reducing   $diam(\pi_k)$   (see Proposition \ref{claim2bis}),  we can obtain  $t_k^j-t_k^{j-1}\le\bar\delta$,  with $\bar\delta$  so small that
\bel{delta1}
\omega_{\tilde m}(LM_f\bar\delta)\le \tilde m(\mu_k)\frac{\varepsilon}{1+2\varepsilon}.
\eeq
Here  $\omega_{\tilde m}$ denotes the modulus of continuity of $\tilde m$,  when restricted  to  $[\mu_{k+1},\sigma]$,  $L$   is  the Lipschitz constant of $U$ on  $U^{-1}([\mu_{k+1},\sigma])$ and $M_f$ is the supremum of $|f|$    on $U^{-1}([\mu_{k+1},\sigma])\times A$.
Hence
\bel{consEk}
\tilde m(U(z(t_k^j))\ge \tilde m(U(z(t))\frac{1+\varepsilon}{1+2\varepsilon},
\eeq
which, together with (\ref{stg}),  implies that
 $$
U(z(t))\le \tilde m^{-1}\left(\frac{2\varepsilon+1}{2\varepsilon+1+t }\,U(x)\right).
$$
 Let us set
\bel{sigma}
\sigma_-(r)\doteq\min\{{\bf d}(z): \ U(z)\ge r\}, \quad
 \sigma^+(r)\doteq\max\{{\bf d}(z): \ U(z)\le r\}.
\eeq
Notice that $\sigma_-$, $\sigma^+:[0,+\infty[\to\R$ are continuous, strictly increasing, unbounded  functions such that $\sigma_-(0)=\sigma^+(0)=0$  and
$$
\forall z\in{U^{-1}([0,\sigma[)}: \quad \sigma_-(U(z))\le {\bf d}(z); \quad  \sigma^+(U(z))\ge {\bf d}(z).
$$
Moreover, it is not restrictive to replace  $\sigma_-(r)$ with $\min\{\sigma_-(r), r\}$.  Let us define  $\beta:[0,+\infty[\times[0,+\infty[\to[0,+\infty[$ by setting
\bel{defbeta}
\beta(r,t)\doteq \sigma^+\circ\tilde m^{-1}\left(
\sigma_-^{-1}(r)
\,\frac{2\varepsilon+1}{2\varepsilon+1+t} \right).
\eeq
Therefore by straightforward calculations it follows that
$$
{\bf d}(z(t ))\le \beta({\bf d}(x), t ) \qquad \forall t\in[0,T_z[.
$$
It  implies  that,  starting from any initial point $x\in U^{-1}(]0,\sigma[)$,
$$
{\bf d}(z(t ))\le \beta({\bf d}(x), t ) \qquad \forall t\in[0,+\infty[.
$$
Let us recall that,  in case $T_z<+\infty$,
we mean that  $z(t)\doteq \bar z$ \, $\forall t\ge T_z$, for some $\bar z \in\partial\C$ (see the footnote to Definition \ref{(GAC)}).
By the arbitrariness of $\sigma>0$, it is easy to extend the construction of $\beta$ from $U^{-1}([0,\sigma])\times[0,+\infty[$   to the whole set $\C^c\times[0,+\infty[$.

\vv
On the one hand, this concludes the proof of part {\bf (i)} of the theorem. On the other hand, the proof of  {\bf (ii)} is  straightforward,  since in view of  (\ref{EP})  we have
$$
 \int_0^{\bar t}l(z(t),\alpha(t))\,dt =\lim_{k\to+\infty}\int_0^{t_k}l(z(t),\alpha(t))\,dt \le \frac{(\varepsilon+1)}{\bar p_0}\, U(x),
$$
which  implies (\ref{Vprop}).

\section{Proofs of some technical results }\label{pprop}

\subsection{Proof of Proposition \ref{claim1}.}    In order to prove (\ref{Hg})  let us observe  that by the definition of $U$,  for any $\sigma>0$,
$U^{-1}([0,\sigma])$ is compact.
Moreover,  for every $\delta\in ]0,\sigma]$  the graph of the restriction of the  set-valued map $x\to D^*U(x)$  to  $U^{-1}({[\delta,\sigma]})$, namely the set
$$
\Gamma_\delta\doteq\left\{(x,p): \   \ x\in U^{-1}({[\delta,\sigma]}), \ p\in D^*U(x)\right\},
$$
is compact. Indeed the set-valued map, $x\to D^*U(x)$ is upper semicontinuous with compact values (see e.g. \cite{AC}).
  Therefore  the continuous function $(x,p)\mapsto H(x,{{\bar p_0}}, p)$  has a maximum on  $\Gamma_\delta$.  For every $\delta\in ]0,\sigma]$, let us set
$$
\hat m(\delta)\doteq -\max_{(x,p)\in\Gamma_\delta} \, H(x,{{\bar p_0}},p)\,\,(>0).
$$
Notice that the function  $\hat m$ is    positive and increasing.  Furthermore, it is  lower semicontinuous.   Finally,   for every $x\in U^{-1}({]0,\sigma[})$    one has
 $$
 H(x,{{\bar p_0}},p)\le -\hat m(U(x)) \qquad\forall p\in D^*U(x).
 $$
The thesis is now proved   by choosing,  for any $\sigma>0$,   a  continuous, strictly increasing,  function $m: [0,+\infty[\to [0,+\infty[$ such that $0<m(r)\leq\hat m(r)$ for every $r\in ]0,\sigma]$  and $m(0)\geq 0$.

 \subsection{Proof of Proposition  \ref{claim2bis}}\label{technical}
Let $(\zeta,\alpha)$ be a trajectory-control pair verifying conditions {\bf (a)--(f)}  of Lemma \ref{claim2} below.   Set $t^0\doteq0$ and,  for every
 $j\in\{1,\dots,\bar n\}$,  define
$$
\tau^j(s) \doteq t^{j-1}+\int_{s^{j-1}}^s \frac{1}{\bar p_0 l(\zeta^j(\sigma),a^j)+ m(U(\zeta^j(\sigma)))}\, d\sigma \quad\forall  s\in[s^{j-1},s^j], \quad t^j\doteq \tau^j(s^j).
$$
Set $\hat t\doteq t^{\bar n}$.
It is trivial to verify that:
\begin{itemize}
\item for every  $j\in\{1,\dots,\bar n\}$,  the path
$$
  z^j: t\mapsto z^j(t) \doteq \zeta^j\circ( \tau^j)^{-1}(t) \qquad t\in[t^{j-1}, t^j].
$$
is a trajectory of the original system in   (\ref{ode})  with initial condition $z^j(t^{j-1})=x^j$,   corresponding to the constant control $a^j$;
\item     the trajectory-control pair $(z,\alpha):[0,\bar t]\to U^{-1}([\hat\mu,\bar\mu])\times A$ given by
$$
(z(t),\alpha(t))\doteq  (z^j(t),a^j) \quad  t\in[t^{j-1}, t^j[ \quad (j\in\{1,\dots,\bar n\}),
$$
  satisfies conditions {\bf (a)--(c)} of Proposition \ref{claim2bis}.
  \end{itemize}

 \begin{Lemma}\label{claim2}
 Let $U$ be a  (MRF), let $\sigma>0$,  and fix a selection $p(x)\in D^*U(x)$.   Let  $m$ be  defined as in Proposition \ref{claim1} when $\sigma$ is replaced with $\sigma+2$, and let $x\mapsto a(x)\in A$ be
 a feedback law\footnotemark  verifying
 \bel{ce}
 \left\langle   p(x)\, ,\, \frac{f(x,a(x))}{\bar p_0\,l(x,a(x))+m(U(x))}\right\rangle\le - 1 \qquad \forall x\in U^{-1}(]0,\sigma+2[).
\eeq
 \footnotetext{Such a feedback exists exactly in view of Proposition \ref{claim1}.}
  Fix $\varepsilon>0$ and   $\bar\mu$, $\hat\mu$ such that   $0<\hat\mu<\bar\mu <\sigma$.

   Then   there exists  $\delta>0$  such that,  for every  partition $\pi=(s^j)$ of $[0,+\infty[$ with diam$(\pi)\le\delta$,    for each $x\in\C^c$ verifying $U(x)=\bar\mu$,   there is a map  $(\zeta,\alpha):[0,\bar s]\to U^{-1}([\hat\mu,\bar\mu])\times A$ verifying
$$
(\zeta(s),\alpha(s))\doteq(\zeta^j(s),a^j) \quad \forall s\in[s^{j-1}, s^j[ \quad (j\ge1),
$$
and a sequence $(x^1,x^2, \dots)\in U^{-1}([0,\bar\mu]) $, where:
 \begin{itemize}
\item[\bf (a)]
$\zeta^1(s^0)=x\doteq x^1$;  for every $j>1$,   $\zeta^{j}(s^{j-1}) =\zeta^{j-1}(s^{j-1})\doteq x^{j}$;
\item[\bf (b)]
for every $j\ge 1$,   $\zeta^j:[s^{j-1},s^j]\to\R^n$ is a solution of the Cauchy problem
$$
\frac{d\zeta}{ds} = \frac{f(\zeta,a^j)}{\bar p_0\,l(\zeta,a^j)+m(U(\zeta))}, \quad
\zeta (s^{j-1}) = x^j,
$$
where
\bel{a}
a^j\doteq a(x^{j});
\eeq
\item[\bf (c)] ${\cal T}^{\hat\mu}_\zeta<+\infty$ and $\bar s\doteq {\cal T}^{\hat\mu}_\zeta$;\footnote{See (\ref{Tz}) for the definition of ${\cal T}^{\hat\mu}_\zeta$. }
\item[\bf (d)] for every $j\ge1$ such that $s^{j-1}<\bar s$,   one has
\bel{dainserire}
U(\zeta^j(s))-U(x^j)\le -\frac{s -s^{j-1} }{\varepsilon+1} \qquad \forall  s\in[s^{j-1}, s^j];
\eeq
\item[\bf (e)]  $U(\zeta(\bar s))=\hat\mu<U(\zeta(s))\le U(\zeta(0))=U(x)=\bar\mu$  \, $\forall s\in[0,\bar s[$,  and
\bel{sceltat2}
\bar s \le (\varepsilon+1)\, U(x).
\eeq
\end{itemize}
Moreover, it is  possible to choose the partition $\pi$ in such a way  that
\begin{itemize}
\item[\bf (f)] $\bar s=s^{\bar n}$ for some integer $\bar n\ge1$,  and,  for every  $ j\in\{1,\dots,\bar n\}$,
\bel{EB}
U(\zeta (s^j))< U(\zeta(s))\le U(\zeta(s^{j-1}) \quad\forall s\in[s^{j-1},s^j[.
\eeq
\end{itemize}
\end{Lemma}

 \noindent{\sc Proof.}
% Let $U$, $\sigma>0$,   $g$,   the feedback law $a=a(z)$, and the constants $\bar\mu$, $\hat\mu\in]0,\sigma[$,    $\bar\mu >\hat\mu$,    be  defined as in the statement.
Fix $\varepsilon>0$ and set
\bel{gdef}
g(x,a)\doteq \bar p_0\, l(x,a)+m(U(x))
\eeq
for all $(x,a)\in U^{-1}({[{\hat\mu/2},{\sigma+1}]})\times A$.  For any continuous function $\phi:{\R^n}\times A\to\R^q$, we use  $M_{\phi}$,   and  $\omega_{\phi}(\cdot)$ to  denote the sup-norm   and the modulus of continuity \footnote{ i.e.,  $$\omega_\phi(r)\doteq\sup\{|\phi(x',a')-\phi(x,a)|: \ \ (x',a'), \, (x,a)\in U^{-1}({[{\hat\mu/2},{\sigma+1}]})\times A, \ |(x',a')-(x,a)|\le r\}.$$}
of $\phi$  in $U^{-1}({[{\hat\mu/2},{\sigma+1}]})\times A$, respectively.
In case   $\phi$ is scalar valued, let  us use  $m_{\phi}$ to denote the minimum of $\phi$ on $U^{-1}({[{\hat\mu/2},{\sigma+1}]})\times A$.
Finally, since $U$ is locally semiconcave,  there exist
$R$, $\rho>0$, $L>0$ such that    for all  $\hat x\in B_n(x,R)\cap U^{-1}({[{\hat\mu/2},{\sigma+1}]})$ one has \footnote{The inequality (\ref{scv}) is usually formulated with the proximal superdifferential  $\partial^P F$ instead of $\partial_C F$. However, this does not make a difference here since $\partial^P F=\partial_C F$ as soon as $F$ is locally semiconcave.}
\bel{scv}
U(\hat x)-U(x)\le  \langle p,\hat x-x \rangle+\rho\,|\hat x-x|^2 \qquad \forall p\in D^*U(x),
\eeq
\bel{Lip}
 |p|\le L \qquad \forall p\in D^*U(x).
\eeq

Let   $\psi:\R^n\to[0,1]$ be a $C^{\infty}$ (cut-off) map such that
\bel{psi}
 \psi = 1 \quad\hbox{on}\quad U^{-1}([{\hat\mu/2}, \sigma]) , \qquad \psi = 0
 \quad\hbox{on}\,\,\, \rr^n\backslash U^{-1}([{\hat\mu/4}, {\sigma+1}])
\eeq

Let us set
$$
\delta\doteq\min\{\delta_1, \delta_2, \hat\mu/2\},
$$
where
\bel{E1}
 \delta_1\doteq \frac{R\,m_g}{M_f},
\eeq
and     $\delta_2>0$   is any positive real number such that
\bel{E2}
L\,\omega_{(\psi\,\frac{ f}{g})} \left(
 \frac{M_f}{m_g} \, r\right) +
\rho\, \frac{M_f^2}{m_g^2} \, r \le \frac{\varepsilon}{\varepsilon+1} \qquad \forall r\in]0,\delta_2].
\eeq

Let $ \pi=(s^j)$ be an arbitrary partition of $[0,+\infty[$ such that  diam$(\pi)\le \delta$.  For each $x\in\C^c$  verifying $U(x)=\bar\mu$, define recursively a sequence of trajectory-control  pairs $(\zeta^j, \alpha^j):[s^{j-1},s^j]\to \R^n\times A$, $j\ge1$,  as follows:
 $$
\zeta^1(s^0)\doteq x^1\doteq x, \ \ a^1\doteq a(x^1);
$$
  for every  $j> 1$,
$$
\zeta^{j }(s^{j-1})\doteq \zeta^{j-1}(s^{j-1})\doteq x^j, \quad a^j\doteq a(x^j);
$$
 for every $j\ge 1$,   $\zeta^j:[s^{j-1},s^j]\to\R^n$ is a solution of the Cauchy problem
$$
\frac{d\zeta}{ds} = \psi(\zeta)\,\frac{f(\zeta,a^j)}{g(\zeta,a^j)}, \quad
\zeta (s^{j-1}) = x^j.
$$
Notice that, by the continuity of the vector field and because of the cut-off factor $\psi$,  any trajectory $\zeta^j(\cdot)$ exists globally  and cannot exit the compact subset $U^{-1}({[{\hat\mu/4},{\sigma+1}]})$.
Let us set
$$
(\zeta(s),\alpha(s))\doteq(\zeta^j(s),a^j) \ \ \forall s\in[s^{j-1},s^j[, \quad \text{for every $j\ge1$.}
$$
 Since, for every $j\ge1$, one has that $s^j-s^{j-1}\le \delta\le \delta_1$,   by (\ref{E1})  it follows that
 $|\zeta^j(s)- x^j|\le R$ \, $\forall s\in [s^{j-1}, s^j]$. Hence, recalling that $|\psi|\le 1$ and  $\psi( x^j)=1$ as soon as $x^j\in U^{-1}([\hat\mu/2,\sigma+1])$,   (\ref{scv}) and (\ref{ce})  imply that, for every  $j\ge1$ such that $s^{j-1}< {\cal T}^{\hat\mu}_\zeta$ (see Definition \ref{Tz}),     one has, $\forall s\in [s^{j-1}, s^j]$,
 $$
\begin{array}{l}
U(\zeta^j(s))-U(x^j)\le
\langle p(x^j),\zeta^j(s)- x^j\rangle+\rho|\zeta^j(s)- x^j|^2 = \\\,\\
\left\langle p( x^j), \int_{s^{j-1}}^{s }{\psi(\zeta^j(\tau)) \,\frac{f(\zeta^j(\tau),a^j)}{g(\zeta^j(\tau),a^j)}\,d\tau} \right\rangle +\rho\left|\int_{s^{j-1}}^{s }\psi(\zeta^j(\tau)) \,\frac{f(\zeta^j(\tau),a^j)}{g(\zeta^j(\tau),a^j)}\,d\tau\right|^2\le   \\\,\\
\left\langle p(x^j),\int_{s^{j-1}}^{s }\left[\psi(\zeta^j(\tau)) \,\frac{f(\zeta^j(\tau),a^j)}{g(\zeta^j(\tau),a^j)}-\frac{f (x^{j},a^j )}{g (x^{j},a^j )}\right]\,d\tau\right\rangle +
\left\langle p(x^j),\frac{f (x^{j},a^j )}{g (x^{j},a^j )} \right\rangle\, (s-s^{j-1})+ \\\,\\ \rho\left( \int_{s^{j-1}}^{s}\left|\psi(\zeta^j(\tau)) \frac{f(\zeta^j(\tau),a^j)}{g(\zeta^j(\tau),a^j)}\right|\,d\tau\right) ^2
\le   \\\,\\
L\,\omega_{(\psi\,\frac{ f}{g})} \left(
 \frac{M_f}{m_g} \, (s^j-s^{j-1})\right)\, (s-s^{j-1})+
\rho\, \frac{M_f^2}{m_g^2} \, (s-s^{j-1})^2-
  (s-s^{j-1})   \le \\ \, \\
\left[L\,\omega_{(\psi\,\frac{ f}{g})} \left(
 \frac{M_f}{m_g} \, (s^j-s^{j-1})\right) +
\rho\, \frac{M_f^2}{m_g^2} \, (s^j-s^{j-1})-1\right] ( s-s^{j-1}).
  \end{array}
$$
Since $ \forall s\in[s^{j-1},s^j]$,  $s-s^{j-1}\le \delta \le \delta_2$,  by (\ref{E2})  it follows that
\bel{stimaE}
U(\zeta^j(s))-U(x^j)\le -\frac{s-s^{j-1}}{\varepsilon+1},
\eeq
which  implies
\bel{stimaE2}
 U(\zeta(s))-U(x)=[U(\zeta^j(s))-U(x^j)]+\dots +[U(\zeta^1(s^1)-U(x)] \le -\frac{s}{\varepsilon+1}.
\eeq
 In particular, (\ref{stimaE2})  yields   that $U(\zeta(s))\le\bar\mu$ for all $s\in[0,s^j]$.

Notice that ${\cal T}^{\hat\mu}_\zeta<+\infty$.  Indeed,
   if   by contradiction  ${\cal T}^{\hat\mu}_\zeta=+\infty$,  (\ref{stimaE2}) held true for all $s\in[0,s^j]$ with   $j$ arbitrarily large, i.e. (since $(s^j)$ is  a partition of $[0,+\infty[$),  for all $s\ge0$. Therefore one would have $\lim_{s\to+\infty} U(\zeta(s))= 0$, which is not allowed, since  \bel{fuori}U(\zeta(s))>\hat\mu\qquad \forall s\in[0,{\cal T}^{\hat\mu}_\zeta[ .\eeq   Let us set
 $$
 \bar s\doteq \bar{\cal T}^{\hat\mu}_\zeta (<+\infty),
 $$
and
 $$
 \quad  \bar n\doteq\sup\{j\ge1: s^{j-1}< \bar s\}.
 $$
 Let us observe that $
 \bar n<+\infty$.

 Finally, notice that, because of (\ref{fuori}), $\psi(\zeta(s) )= 1$ for every $s\in[0, s^{\bar n}]$. Hence,  for any $j\in\{1,\dots,\bar n\}$,  $\zeta^j(\cdot)$ is  a solution of
 $$
 \frac{d\zeta}{ds} = \frac{f(\zeta,a^j)}{g(\zeta,a^j)} \ \ \forall s\in[s^{j-1},s^j], \quad
\zeta (s^{j-1}) = x^j.
$$
It follows that conditions  {\bf (a)--(e)}  are  satisfied.  Notice however that in general {\bf (f)} does not hold. Indeed it may happen that $s^{\bar n}>\bar s$. In addition, the first inequality of (\ref{EB}), namely
$U(\zeta (s^j))< U(\zeta(s))$, may fail to be verified for some $s\in]s^{j-1},s^j[$.

In order to prove {\bf (f)}, it is sufficient to slightly refine the previous construction:
\begin{itemize}
\item  In case  (\ref{EB}) does not hold in $[0,\delta]$, redefine $s^1$ by setting $s^1\doteq\inf\{s\in]\hat 0,\delta]: \  U(\zeta^1(s))\le U( \zeta^1(\delta))\}$.
\item  For every $j>1$, choose $a^j$ and  $\zeta^j:[s^{j-1},s^{j-1}+\delta]\to\R^n$ with the same procedure we have followed in the previous construction. In case  (\ref{EB}) does not hold in $[s^{j-1}, s^{j-1}+\delta]$, set $s^j\doteq\inf\{s\in]s^{j-1},  s^{j-1}+\delta]: \  U(\zeta^j(s))\le U(\zeta^j(s^{j-1}+\delta))\}$.
\end{itemize}

The trajectory  $\zeta$ defined by setting  $\zeta(s)\doteq \zeta^j(s)$ \, $\forall s\in[s^{j-1},s^j[$, $j\ge1$,  verifies (\ref{EB}). It remains to prove that $(\bar s=)\, {\cal T}^{\hat\mu}_{\zeta}=s^{\bar n}$ for some integer $\bar n$. Begin with observing that $U(x^{j+1})=U(\zeta^j(s^j))=
U(\zeta^j(s^{j-1}+\delta))$. Hence, by (\ref{stimaE}) it follows that
$$
U(x^{j+1})-U(x^j)= U(\zeta^{j }(s^{ j-1}+\delta))-U(\zeta^{j }(\ s^{j-1}))\le -\frac{\delta}{\varepsilon+1},
$$
so that
$$
U( x^{j+1})\le U(x) -\frac{j\delta}{\varepsilon+1}.
$$
Since $j\in\N$,  we get $ s^{\bar n}={\cal T}^{\hat\mu}_{\zeta}$ for some integer   $\bar n$ smaller than the first $n$ such that $n\delta> (\varepsilon+1)(\bar\mu-\hat\mu)$.

%%%%%%%%%%%%%%%
\section{ A remark on supersolutions }\label{Concl}
The notion of  (MRF)  can be restated by replacing the strict  inequality (\ref{C1}) with a supersolution condition.

Preliminarly, let us  recall some basic facts from nonsmooth analysis.  We remind that we are using $\partial_{{\bf C}}F$ and $D^*F$ to denote the Clarke's generalized gradient and the set of limiting gradients, respectively (see Definition \ref{D*}).

\begin{Definition}\label{D-}
Let $\Omega\subset\R^n$ be an open set, and let  $F:\Omega\to\R$ be a locally bounded  function.   For every $x\in\Omega$,  the set
$$
 \quad D^-F(x)\doteq\left\{p\in\R^n: \ \ \liminf_{y\to x}\frac{F(y)-F(x)-\langle p,(y-x)\rangle}{|y-x|}\ge0\right\},
$$
is called the {\it subdifferential of $F$.}
\end{Definition}

We recall that  $D^-F(x)$ is a  closed, convex (possibly empty) set. If $F$ is differentiable at $x$, then  $D^-F(x)=\{\nabla F(x)\}$.  Moreover, when $F$ is locally Lipschitz,  $D^-F(x)\subset \partial_{{\bf C}}F(x)=co\, D^*F(x)$.

 \begin{Proposition}\label{L1} Let $U:\overline{\C^c}\to\R$ be a (MRF) and let $\bar p_0\ge 0$ be the constant for which  (\ref{C1}) holds true.  Then the strict  inequality (\ref{C1}) can be equivalently replaced by the following  condition:
 \begin{itemize}
 \item[]
 for every $\sigma>0$,  there exists a continuous, strictly increasing   function  $m=m_\sigma:[0,+\infty[\to\R$ verifying $m(r)>0$ \, $\forall r>0$,  such that   $U$  {\em is a  viscosity supersolution of equation $-H(x,{{\bar p_0}}, Du)- m(u)=0$
 in $U^{-1}(]0,\sigma[)$}, namely,    one has
\bel{ssol}
 H(x,{{\bar p_0}}, D^- U(x))\le- m(U(x)) \quad \forall x\in U^{-1}(]0,\sigma[).
\eeq
\end{itemize}
 \end{Proposition}
\noindent {\sc Proof.} \ In view of Proposition \ref{claim1},  in order to show that  (\ref{ssol}) implies (\ref{C1})   it is enough to prove that,  for any $\sigma>0$,  (\ref{ssol}) implies
\bel{D*sol}
 H(x,{{\bar p_0}}, D^*U(x))\le- m(U(x)) \quad \forall x\in U^{-1}(]0,\sigma[)
\eeq
(for the same function $m$ as in  (\ref{ssol})).  For any $x\in U^{-1}(]0,\sigma[)$ and $p\in D^*U(x)$,   there is a sequence $(x_n)\subset U^{-1}(]0,\sigma[)\cap\, DIFF(U)$ such that \linebreak $\lim_n(x_n, \nabla U(x_n))=(x,p)$.   Since
\bel{D-}
 D^-U(x)=\{\nabla U(x)\} \quad\forall x\in DIFF(U),
  \eeq
  by hypothesis (\ref{ssol})  one has
$$
H\left(x_n,\,{{\bar p_0}},\, \nabla U(x_n)\right)\le -m(U(x_n)),
$$
 for each  natural number $n$.  Passing to the limit  an $n$ tends to infinity  we get (\ref{D*sol}).  The converse implication  is a straightforward consequence of the following  relations (see e.g. \cite{CS}):
\bel{D-}
\begin{array}{l}
D^*U(x)=D^-U(x)=\{\nabla U(x)\} \quad\forall x\in DIFF(U);  \\\,\\
D^-U(x)=\emptyset  \quad\forall x\in \C^c\setminus DIFF(U).
\end{array}
\eeq

\vv
 {\em Acknowledgments.} The authors  thank Fabio Priuli and the anonymous referees for carefully reading the paper and making useful comments.

 \end{document}